\documentclass[fleqn,10pt]{wlscirep}
\usepackage[utf8]{inputenc}
\usepackage[T1]{fontenc}
\usepackage{float}
\title{Simulation Models for Sustainable, Resilient, and Optimized Global Electric Vehicles Supply Chain}

\author[1]{Tareq Alsaleh}
\author[1,*]{Bilal Farooq}
\affil[1]{Laboratory of Innovation in Transportation (LiTrans), Toronto Metropolitan University, Toronto, M5B 2K3, Canada}
\affil[*]{Bilal.Farooq@torontomu.ca}


\begin{abstract}
 While the transition to electric vehicles (EVs) is essential for decarbonizing the transportation system, the production and distribution of EVs entail substantial carbon costs. To ensure these emissions are accurately accounted for and effectively mitigated, this research introduces a digital twin of the EV's supply chain, addressing a critical gap in current EV life cycle analyses and providing the first comprehensive quantification of its environmental sustainability and resilience. This simulation model replicates global market dynamics and captures the complexity and uncertainty of the EV supply chain, enabling a thorough evaluation of its carbon footprint, sustainability, resilience, and what-if counterfactual scenarios for alternative market structures. The results reveal that average supply chain emissions range from 6.42 to 6.94 Kg e-CO$_2$/KWh across different battery technologies. Additionally, the mass flow analysis shows unbalanced dependencies at all supply phases, with one geographical region significantly dominating the supply chain structure, highlighting the current supply chain architecture's low resilience and high vulnerability. In light of these findings, the study introduces an optimization model for hub and resource allocation configuration, effectively reducing vulnerability levels and supply chain emissions by up to 80\%.
\end{abstract}

\begin{document}

\flushbottom
\maketitle
%
%
\thispagestyle{empty}

\section*{Introduction}\label{Sec1}

As the electrification of passenger vehicles became the central focus in the transition to zero-emission transportation systems, governments worldwide have pledged to achieve 100\% electrification by the year 2050 or earlier  \cite{IEA2023EV}, which in turn led to a significant surge in the global sales of electric vehicles (EVs) \cite{IEA2023ElectricCarStock}. 
 This shift towards EVs is not a standalone strategy but rather part of a broader initiative geared toward adopting environmentally sustainable practices across all industries. Central to this ambitious agenda, which ultimately aims to combat climate change and regulate the pace of global warming, are two fundamental considerations: the effectiveness of transportation electrification in curbing the pace of global warming and the establishment of a secure and sustainable supply chain necessary to meet the manufacturing and distribution targets of EVs. 

The effectiveness of EVs as an environmentally sustainable option to replace internal combustion engine vehicles (ICEVs) is often addressed in literature through comparative life cycle assessments \cite{nimesh2021implication, koroma2020prospective, qiao2020life, almeida2019quest,kim2019environmental, shen2019china, petrauskiene2020comparative, liu2020optimal, petrauskiene2021comparative, cusenza2019energy, wang2022review}. Among these studies, estimated emissions from the EV's life cycle ranged from 100 to 253 g e-CO$_2$/km and from 170 to 340 g e-CO$_2$/km for the ICEV \cite{qiao2020life, kim2019environmental}. Variations in emissions are attributed to differences in system boundaries, study vehicle types, functional units, vehicles, battery lifetime, carbon intensity of power generation, and driving environment \cite{xia2022life}.

Most of the studies evaluated EVs in the Chinese \cite{shi2019life, qiao2017comparative},  United States (US) \cite{samaras2008life, challa2022well}, and European \cite{pipitone2021life, tagliaferri2016life} contexts. This location-specific analysis is critical to capture the variation in different power mix scenarios for charging EVs, where the percentage of renewable energy plays a critical role in reducing EV life cycle emissions \cite{del2018life}. For instance, to demonstrate the power mix scenarios role in determining the effectiveness of the electrification strategy, a well-to-wheel (WTW) study in China showed that EVs reduce emissions by 10 to 20\% in northern provinces compared to ICEVs, while the reduction could amount to 50\% in southern provinces  \cite{shi2019life}. As for the US, Samaras et al. \cite{samaras2008life} concluded that plug-in hybrid vehicle PHEV reduces 32\% emissions compared to ICEV using a cradle-to-grave life cycle emission analysis. Another study in the US context concluded that ICEVs emit higher emissions than EVs, based on the WTW analysis for the annual travel distance \cite{challa2022well}. 

EV's environmental impact studies also vary in terms of system boundaries. Several studies analyze cradle-to-grave emissions \cite{shi2019life, samaras2008life, pipitone2021life, tagliaferri2016life}, while others only examine the production phase known as cradle-to-gate analysis. For instance, Qiao et al. \cite{qiao2017comparative}, utilizing the literature on energy consumption, reported that the emissions from the EV production range between 14.6 and 14.7 tonnes of e-CO$_2$, around 60\% higher than the ICEV.  
Battery and attachment emissions accounted for 2.788 and 2.892 tonnes of e-CO$_2$ per Nickel-Manganese-Cobalt (NMC) and Lithium-Iron-Phosphate (LFP) batteries respectively, comprising the highest percentage of total emissions after the vehicle body (exterior and interior components combined). A Fixed battery type and size were used in the study, considering a battery mass of 170 and 230 Kg per NMC and LFP batteries. The battery assembly emissions were assumed to be constant at 141.5 Kg e-CO$_2$ for both types. Using fixed battery type and mass is common in evaluating EV emissions at the production phase. Another study considered NMC and LFP batteries with battery masses of 214 and 273 Kg, respectively, and concluded that battery production accounts for 35 to 41\% of the global warming potential emissions of the EV production phase \cite{hawkins2013comparative}. 

In the process of evaluating EV emissions, the supply chain and assembly are difficult to assess in detail due to uncertainties in the supply links and the unavailability of accurate data. Previous studies used linear functions of mass and scalar factors to estimate energy consumption and emissions \cite{ma2012new} or assign fixed emission values to the supply chain and assembly process \cite{qiao2017comparative}. The supply chain is not commonly considered in the evaluation due to the unclarity of data, shipping situation and transport distances \cite{wang2021comparative}. In addition, the supply chain is often considered a minor contributor to emissions in environmental assessment \cite {wu2021environmental}, which could be the probable reason why most studies prefer to only mention the transport links \cite{ballinger2019vulnerability}. 

The consideration and analysis of supply chain emissions are also missing from the literature investigating the sustainability of the electrification strategy application, especially in terms of manufacturing scale-up, resource distribution and the allocation of manufacturing hubs for building a secure and resilient supply chain architecture. The sustainability of implementing the strategy is often only sought from the analysis of critical mineral reserves to meet the global demand for EVs and other industrial uses \cite{jones2020ev} without serious regard for the emissions of the supply chain or the potential for optimization to reduce cost, travel distance and time, and build a resilient global supply chain that is less vulnerable to disruption and instabilities.

A large body of literature investigated the production and operation phases to assess EV emissions \cite{nimesh2021implication, koroma2020prospective, qiao2020life, almeida2019quest,kim2019environmental, shen2019china, petrauskiene2020comparative, liu2020optimal, petrauskiene2021comparative, cusenza2019energy, wang2022review, shi2019life, samaras2008life, pipitone2021life, tagliaferri2016life, hawkins2013comparative, qiao2017comparative, ma2012new, wang2021comparative, wu2021environmental, jones2020ev,xia2022life,challa2022well,ballinger2019vulnerability,del2018life}. Evaluations often consider a specific type and size of vehicle and battery sets, mainly to obtain the energy consumption rates for the production of different parts of that specific model. Conversely, due to the uncertainty of the data and shipping scenarios for raw materials and finished products, the emissions from the supply chain are regarded as too difficult to obtain, and therefore assembly and supply chain emissions are often added as a fixed value or as a linear function of mass regardless of the supply chain scenarios. This simplification reveals a substantial gap in our understanding of how present and future supply chain structures and market dynamics influence the sustainability of EV adoption. Furthermore, there is a notable lack of insight into the overall efficiency and resilience of the existing supply chain. An in-depth understanding of the supply chain and mass flow characteristics is essential not only for sustainability analysis but also to identify vulnerabilities and optimize efficiency. Resilience analysis provides insights into the ability of the supply chain to withstand geopolitical tensions, natural disasters, and market fluctuations. 

 
This research addresses the two fundamental questions raised earlier by providing a detailed EV supply chain analysis through simulations of the EV market structure. The resulting digital twin is used to explore and test various what-if market scenarios to improve the resilience and sustainability of the EV supply chain and optimize manufacturing hub locations and their sourcing options. Ultimately, supply chain resilience analysis aims to establish a market structure with an uninterrupted supply chain by providing sustainable expansion options to boost global production capacity and meet electrification targets.
Additionally, to help place the significance of the supply chain emissions compared to the overall pre-operation emissions, we provide a generalized method to evaluate EV battery emissions from the mining of the raw materials as a function of battery capacity in kWh, utilizing the literature on individual mineral analysis, without reference to specific car size or model as used in previous studies.

\section*{Results}
This section presents the simulation model outcomes on the supply chain's sustainability analysis, considering various battery technologies and detailing the emissions across different transport modalities via land and sea shipments. We also analyze emissions data on a per-company and per-market basis, providing insights into the impact of the supply structure on environmental sustainability and resilience. Additionally, we present the mass flow data for detailed resilience and vulnerabilities analysis and potential improvement areas. Furthermore, this section discusses the strategic allocation of resources and the positioning of manufacturing hubs, demonstrating how resource management can significantly reduce emissions and bolster overall supply chain resilience.

\subsection*{Global Supply Chain Sustainability Analysis}\label{subsec2.1}

\begin{figure}
  \centering
  \includegraphics[height= 15cm, width= 18cm]{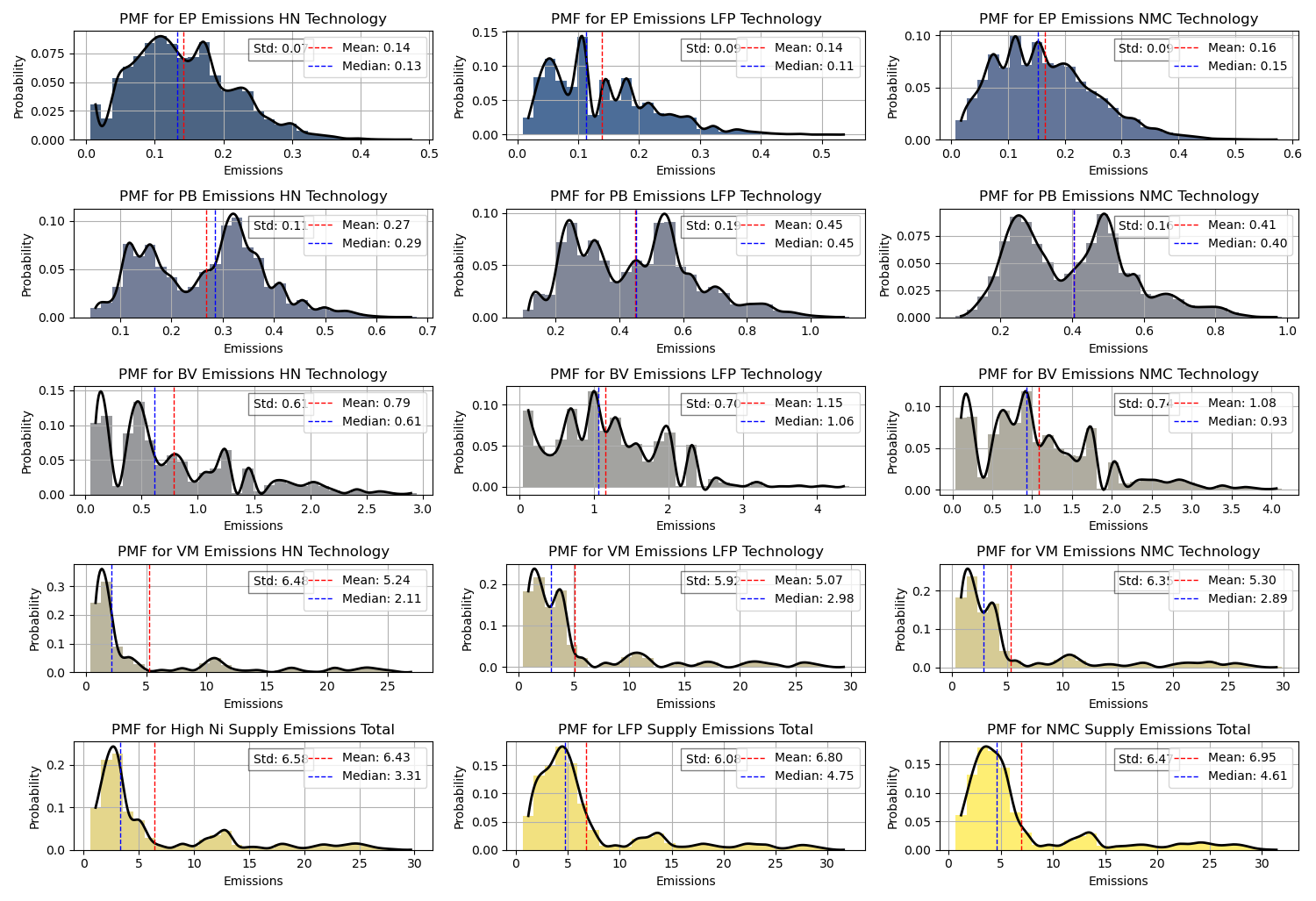}
  \caption{Probability Mass Function of Supply Chain Emissions for Different Supply Phases for Different Battery Types}
  \label{fig:1}
\end{figure}

The Supply Chain for EV production Simulator (SCEV-Sim) was developed in order to track, record and analyze the mass flow for EVs production; the details on SCEV-Sim formulation and the simulation input are provided in the Methods Section. This simulation model serves as the foundation of our investigation into the sustainability of the EV supply chain, considering the e-CO$_2$ emissions—a key metric for assessing environmental impact and sustainability.
;
The sustainability analysis results are segmented into four supply phases: extraction to processing (EP), processing to battery production (PB), battery production to vehicle production (BV) and vehicle production to market distribution (VM). For each phase, we present the emission values through it's probability mass function (PMF) per battery technology, as shown in Figure \ref{fig:1}. The SCEV-Sim reveals that the average EV supply chain emissions per KWh vary between 2 and 8\% depending on the battery technology used in the production. Vehicles with NMC battery technology exhibit the highest mean emissions, averaging 6.87 kg e-CO$_2$ / KWh. This is followed by LFP and high nickel batteries, with average emissions of 6.71 and 6.37 kg e-CO$_2$ / KWh, respectively. The PMFs also show that emissions could peak at 30 kg e-CO$_2$ / KWh, albeit with a substantially lower frequency. 

Our findings indicate that global supply chain emissions markedly contribute to the overall pre-operation emissions for EV production, contrary to prevalent assumptions. This revelation highlights the importance of integrating a detailed supply chain analysis within the frameworks of EV life cycle assessments. Additionally, the results call for an in-depth investigation to dissect the emissions across different transport modes, suppliers, and markets. Such analyses are essential to inform the optimization scenarios and strategies that ensure the sustainability of the electrification approach, especially with the global rise in demand for EVs. 

\begin{figure}
  \centering
  \includegraphics[width=\textwidth]{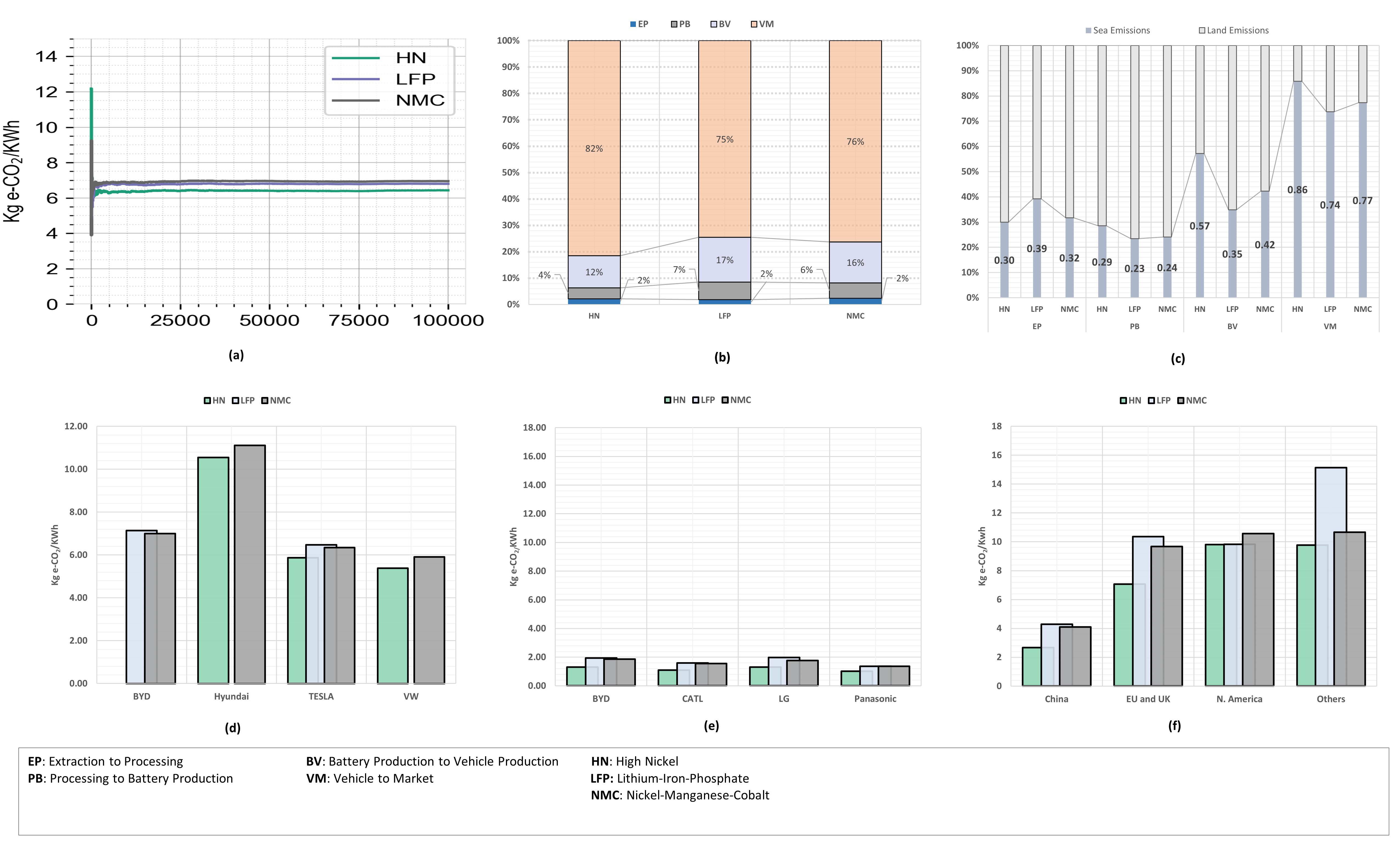}
  \caption{ EV Supply Chain Analysis Emissions Breakdown \textbf{(a)} Cumulative Average per Supply Chain Emission per KWh per Battery Type \textbf{(b)} Percentage Emission per Supply Phase \textbf{(c)} Percentage Supply Emission per Transport Mode \textbf{(d)} Average Supply Emission per KWh per Car Manufacturer \textbf{(e)} Average Supply Emission per KWh per Battery Manufacturer \textbf{(f)} Average Supply Emission per KWh per Market}
  \label{fig:2}
\end{figure}

For a more granular understanding of the emissions, Fig. \ref{fig:2} provides a breakdown of the emissions by phase, transport modality, manufacturers of cars and batteries, and markets, as well as the simulation stabilization results.  Fig. \ref{fig:2} (a) shows the evolution of the cumulative average of the supply chain for vehicles with the three variants of battery types achieving stable cumulative average results at around 10K iterations. However, to further investigate the asymptotic change in the cumulative average, 100K and one million iterations were tested. The change in cumulative averages between 10K and 100K was around 1\%, while the change between 100k and one million was between 0.1\% and 0.3\%. These recorded changes were not necessarily in the same direction, which precludes the possibility of an asymptotically decreasing or increasing pattern and reaffirms the stability of the model output. Detailed results for 10K, 100K, and one million iterations are presented in the Supplementary Information File Fig. S1, S2 and S3. 

Fig. \ref{fig:2} (b) shows that a substantial portion of the emissions are caused by the final phase of the supply chain, which accounts for 75\% to 82\% of total emissions. The percentage contribution of each phase gradually increases as we move from raw minerals to finished products, corresponding to the increase in mass flow. 
Fig. \ref{fig:2} (c) disaggregates the emissions by transport type, where the significant emission contribution from land transport can be noted for battery production phases (EP and PB) before the sea transport takes over as a leading emitter for the delivery phase to markets. 

In the context of the car manufacturers shown in Fig. \ref{fig:2} (d), there is a clear disparity among companies, with Hyundai exhibiting notably higher average supply emissions per kWh than its counterparts. Such disparity is not as clear on the battery manufacturers' side as shown in Fig. \ref{fig:2} (e) of the same figure. However, on average, Panasonic recorded around 20\% less in supply chain emissions than other prominent battery manufacturers. Fig. \ref{fig:2} (f) displays the breakdown of the supply chain emissions per the largest consumer markets. The results show that the Chinese market had considerably lower average emissions from the supply chain phase, followed by the EU for the supply of vehicles with high Nickel and NMC batteries and the US for the vehicles with LFP batteries. It is important to note here that the values are presented per kWh and do not resemble the total emissions from each market but rather an average for each KWh supplied to the market, regardless of the total quantities.     
   
\begin{figure} [H]
  \centering
  \includegraphics[width=\textwidth]{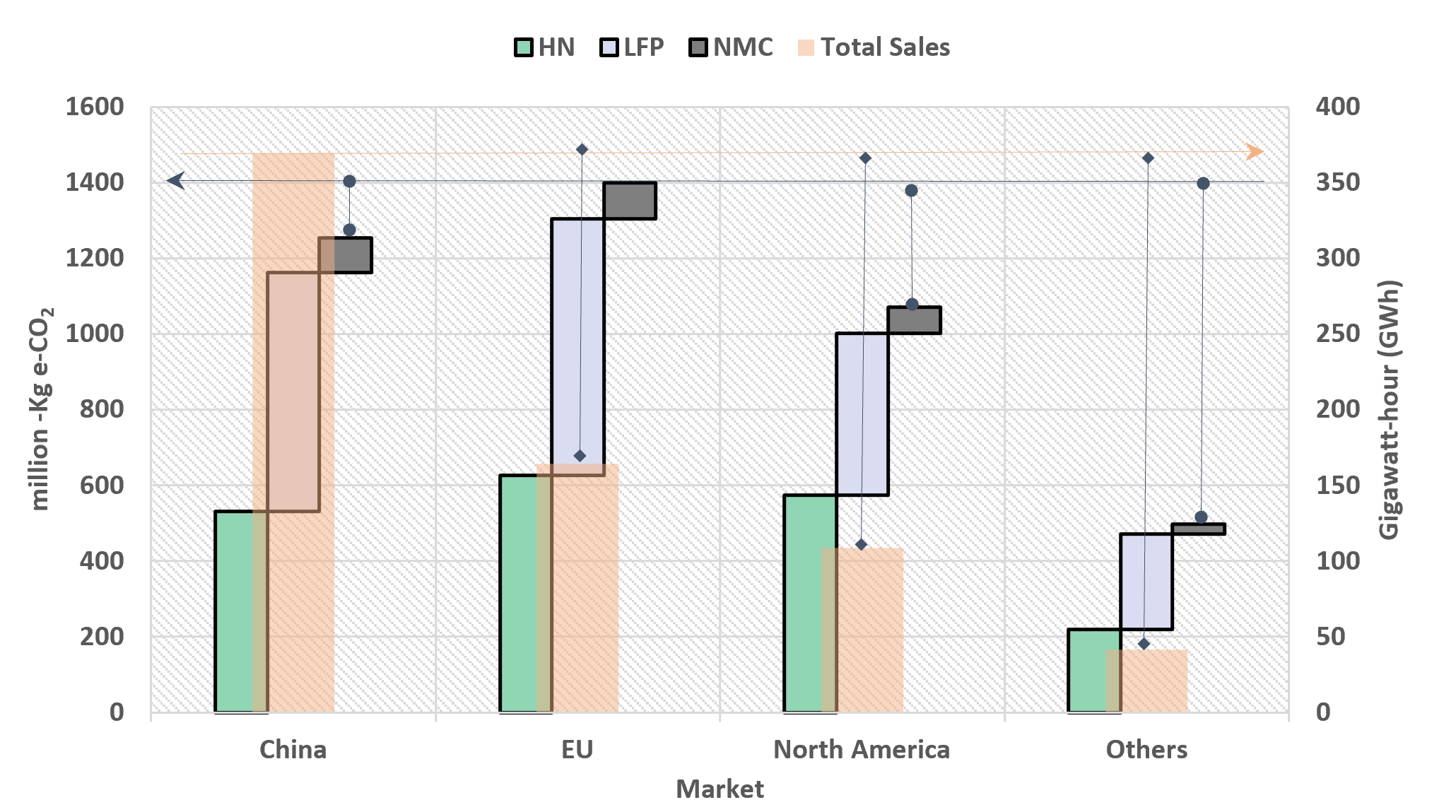}
  \caption{ Total Emissions Per Market Per Battery Technology (Y-Axis Left)- Total Electric Vehicles Sales per Market in GWh (Y-Axis Right) }
  \label{fig:3}
\end{figure}

As for the total emissions, Fig \ref{fig:3} consolidates the emissions into an aggregate level, providing an overview per market per battery technology, considering the recorded global sales for 2023. The left Y-axis of the figure presents the total emission values per market in a million kgs of e-CO$_2$, while the right Y-axis show the total sales of electric vehicles in GWh. While the EV sales in China are 220\% and 340\% higher than the EU and US, respectively, China's total EV emissions from the supply chain are 11\% lower than the EU and only 17\% higher than the US. This demonstrates the significant role of China's localized supply chain and the impact of its minimal average supply chain emission per kWh shown in Fig. \ref{fig:2} (f) on the overall sustainability of its EV supply chain. 

\subsection*{Global Supply Chain Mass Flow and Resilience Analysis}\label{subsec2.2}

\begin{figure}
  \centering
  \rotatebox{270}{%
    \begin{minipage}{\textheight}
      \centering
      \includegraphics[width=\textwidth]{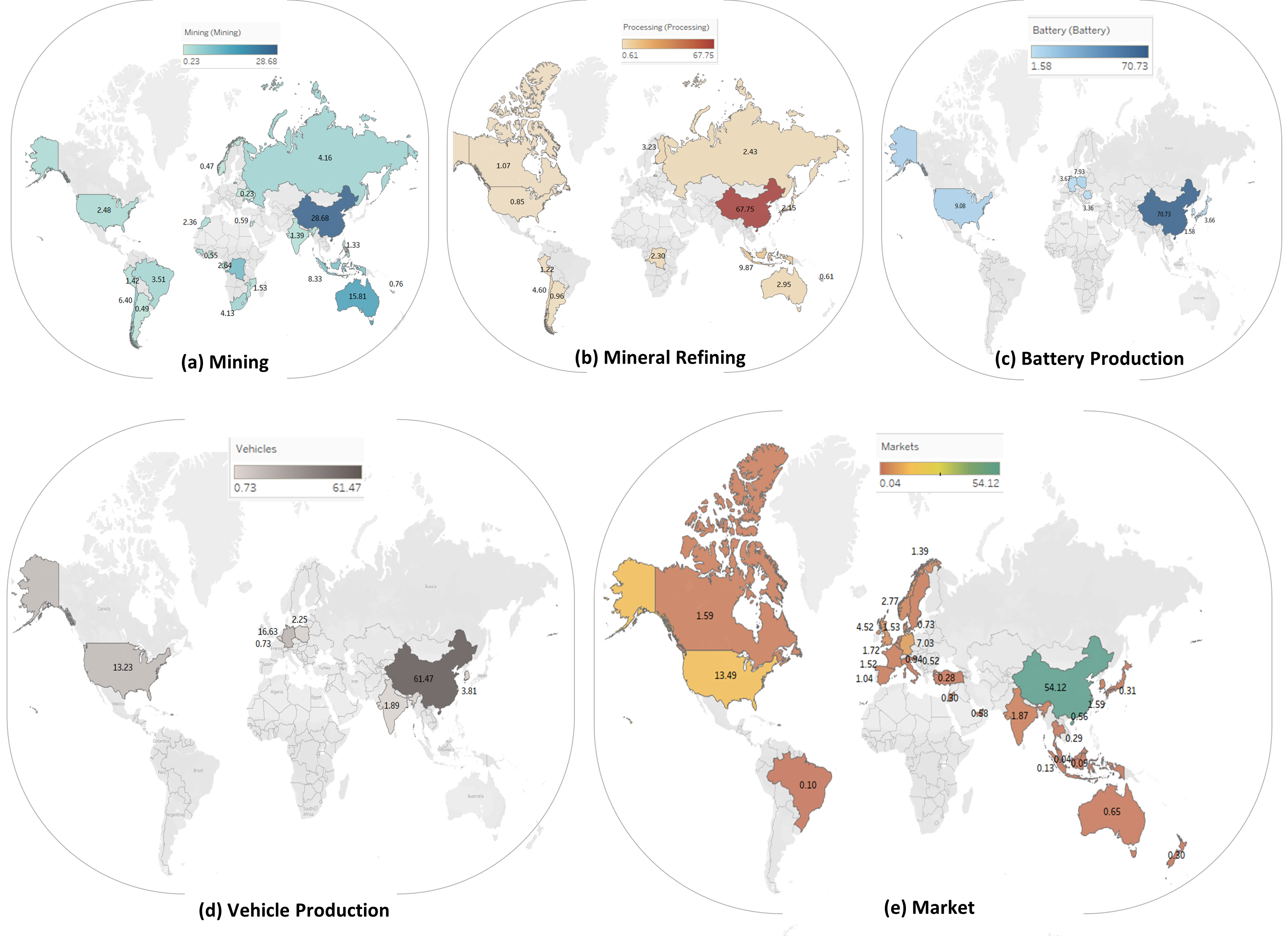}
      \caption{Supply Chain Resilience Analysis- Market Share per Country Per Supply Phase \textbf{(a)} Mining Phase  \textbf{ (b)} Mineral Refining Phase \textbf{(c)} Battery Production Phase \textbf{(d)} Vehicle Production Phase \textbf{(e)} Market}
      \label{fig:4}
    \end{minipage}%
}
\end{figure}
\begin{figure}
  \centering
  \rotatebox{270}{%
    \begin{minipage}{\textheight}
      \centering
      \includegraphics[width=\textwidth]{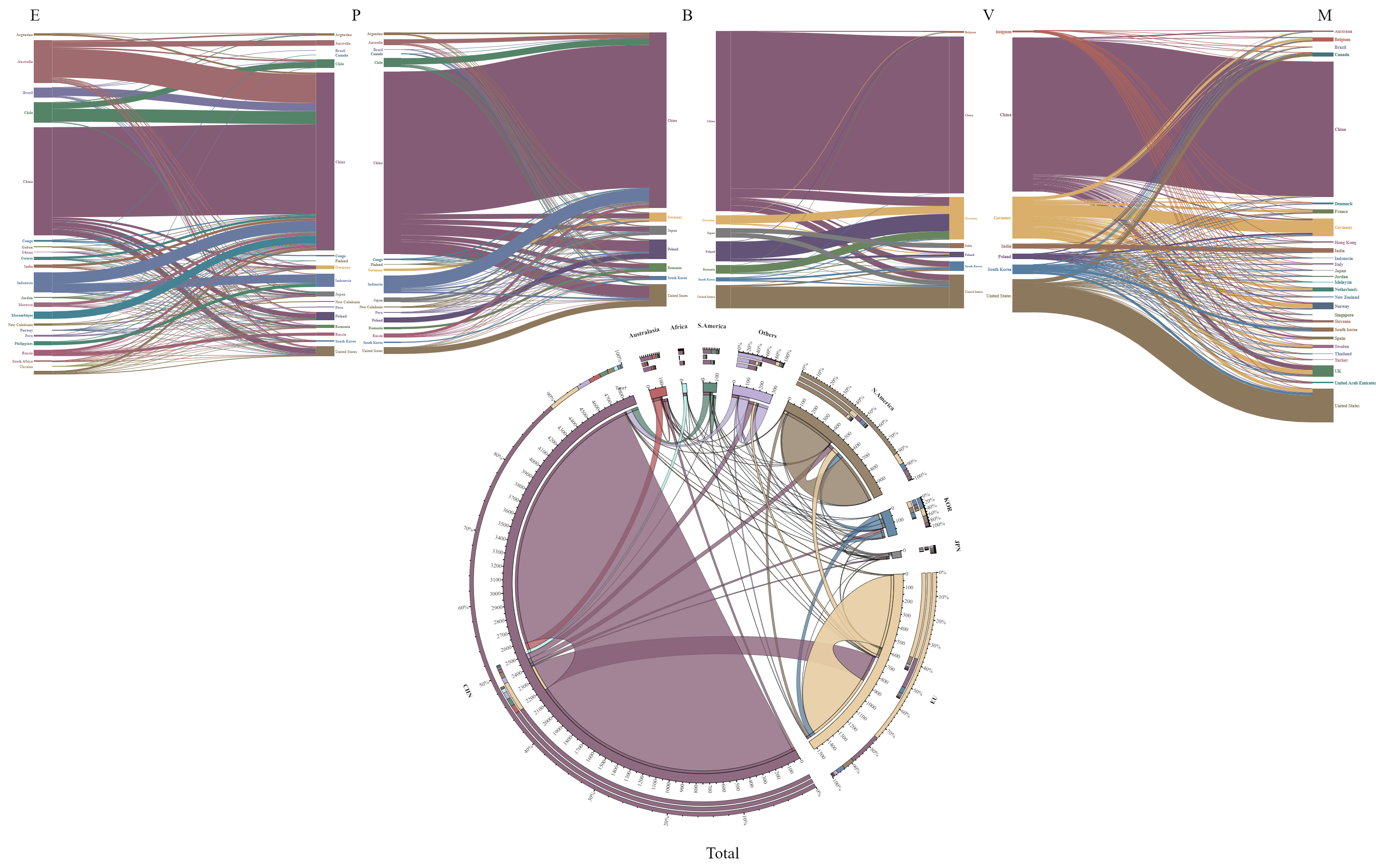} 
      \caption{Mass Flow Tracking and Analysis. \textbf{E:} Mineral Extraction to \textbf{P:} Mineral Refining and Processing to \textbf{B:} Battery Production Hubs to \textbf{V:} Vehicle Production Hubs to \textbf{M} Markets and the Total Mass Flo+w Transfer for the Full Supply Cycle in Kg/ 100 KWh.}
      \label{fig:5}
    \end{minipage}%
}
\end{figure}

SCEV-Sim was also used to track and quantify the supply chain's resilience. The resilience was examined using two indicators: 1) The market share of production volumes per source, which quantifies the production capacities and significance at different supply chain phases. 2) The origin-destination (O-D) mass flow in Kg transferred per KWh, identifying production, consumption, and the supply links used for the mass flow between supply phases. The results from both indicators are presented in Fig. \ref{fig:4} and Fig. \ref{fig:5}.

Fig. \ref{fig:4} (a) shows the mining production per source, including all minerals required for EV production considering high nickel, LFP and NMC battery variants. Fig. \ref{fig:4} (b), on the other side, shows the geographical locations for mineral refining to the grade usable for EV battery production. Panel (c) and (d) of Fig \ref{fig:4} map out the major production hubs for battery and vehicle, respectively, according to their production volumes compared to the total global production. Finally, Fig. \ref{fig:4} (e) presents the market share of EV sales.

The highest diversity in the market structure distribution is at the mining stage. Worth noting that this is based on the overall minerals required for EV production and not on minerals classified as critical. However, there is clear dominance for China for the remainder of the supply phases, with a share that amounts to 70.7\% of the origins of EV batteries in the battery production phase and a home of 67.7\% of the refining market, while it has a slightly lower but still dominant market share on the vehicle production side with 61.4\%. The graphs further demonstrate that China is the only location without a deficit between consecutive production stages, highlighting its supply chain resilience and capacity to satisfy the entire supply chain cycle locally—except for mining. 

As for the two other major EV markets, the US and the EU, the US is present for all stages to varying degrees but with an actual competitive advantage at the vehicle production stage. On the other hand, the EU is almost absent from two critical stages, namely mining and refining, with a humble presence on the eastern European side. However, the EU's contribution to the EV supply chain picks up on the battery and vehicle production sides, where it is home to approximately 15\% and 19\% of produced batteries and vehicles, respectively. However, the analysis reveals the limited resilience of the EU and the US EV production markets, particularly in their capacities for mineral refining and battery production, which lag behind their vehicle production capacities and market demands. This discrepancy affects their ability to build a robust and resilient supply chain that meets their local demand and competes effectively in the global EV market.  

To further assess the resilience of the global EV supply chain, the O-D mass flow was recorded using SCEV-Sim and presented in Fig. \ref{fig:5}. The analysis shows that 4103 kg of mass is transferred for every 100 KWh of EV production, with 80\% of the mass transfer occurring domestically or within the same region. In terms of share distribution, China dominates the global mass flow, taking part in 65\% of the global flow, with 80\% of this movement taking place domestically. North America and the EU contribute 14.5\% and 23\% of the mass flow, respectively, with 67\% and 60\% domestic mass flows in these regions. 

While the percentage of the mass flow passing through each region is indicative of the region's significance on the EV production landscape, and the percentage of the domestic mass transfer reflects the regional supply chain's self-sufficiency, neither metric necessarily provides a comprehensive evaluation of the overall supply chain resilience for different regions. To effectively assess the resilience of supply chains for different regions, their domestic and international mass flow needs to be evaluated in tandem with the flow balance. For instance, China exhibits a positive flow balance of 140 Kg per 100 KWh, compared to a negative flow balance of 102 and 190 Kg per 100 KWh for both North America and the EU, respectively. The flow balance here provides insights on whether the percentage of the non-domestic flow complements the supply chain resilience or, in contrast, reveals its vulnerability. 

\subsection*{Global Supply Chain Optimization and Resource Allocation}\label{subsec2.3}

\begin{figure}[ht]
\centering
\includegraphics[width=1\linewidth]{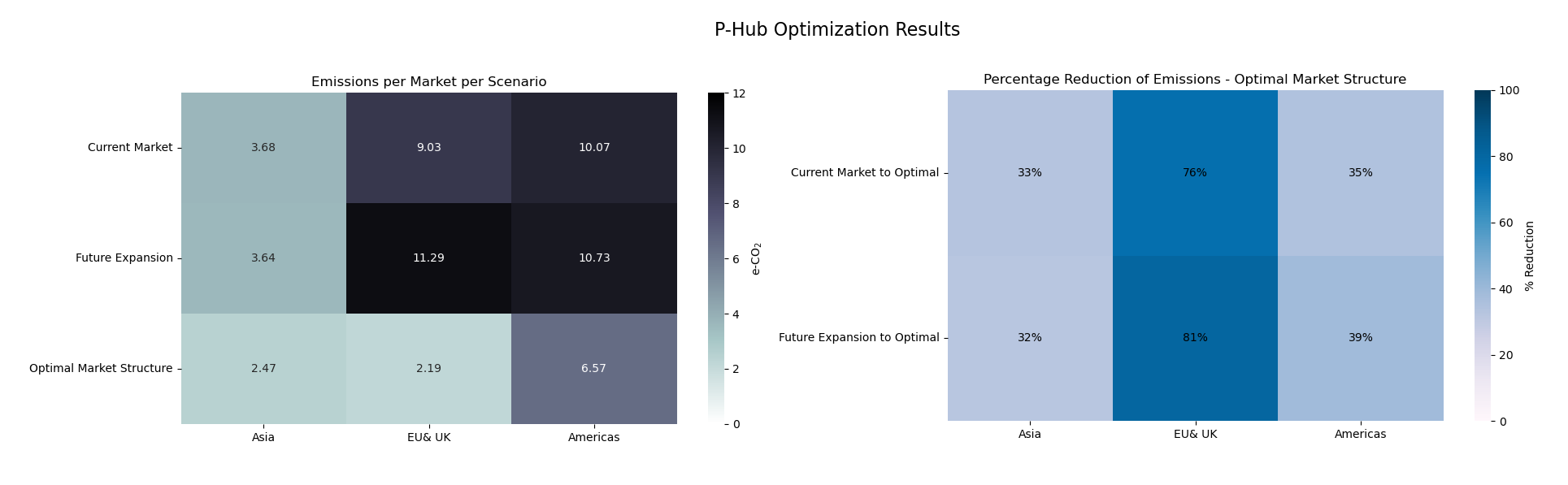} 
\caption{P-Hub Optimization Results}
\label{fig:10}
\end{figure}

In light of the sustainability and resilience analysis findings, there is an evident need to optimize the supply chain and the market structure to reduce emissions and promote a sustainable expansion of production capacities. In this study, a P-hub optimization was formulated to define production hub locations while also allocating appropriate resources necessary for the EV production, as detailed in the Methods sections. The results of the optimization problem are presented in Fig. \ref{fig:10}. As can be noted, the optimization of the hubs at two production levels and the resource allocation result in a significant reduction of emissions. The emissions from the optimized supply structure were compared with both current market structure as well as the near future anticipated market structure emissions, where increased international trade levels are expected. The results show that the optimized supply chain has the lowest average emission of 2.19 kg e-CO$_2$ / KWh for the EU market, followed by the Asian and the American markets with average emission values of 2.47 and 6.57 kg e-CO$_2$ / KWh, respectively. While, the EU market exhibit the lowest average emissions in an optimized supply structure, it certainly has the highest average emissions for the anticipated future market structure with average emission value of 11.29  kg e-CO$_2$ / KWh a 25\% expected increase from its current average emission values. 

The analysis results from the anticipated global expansion show an increase in emissions per KWh, which is expected as new markets and hubs emerge. The increase is more substantial for the EU market as result of an increase in the international trade scenarios to meet their domestic needs, especially the increased trade the Chinese EV manufacturers. 

The optimization shows a high percentage reduction in the emissions for the EU markets and a lower but still significant decrease for the Americas markets; this is basically a result of the vast distance covered by the Americas market, as it serves both North America and South America. Considering separate production hubs with their associated optimal resource allocation for these markets is expected to yield a further reduction in emissions. However, the structure is maintained as is in the current analysis to adhere to historical trends of mass production capacities and to provide viable options for global expansion of production hubs at different levels. 

\section*{Discussion}

The supply chain emissions values are of significance and must be considered in EV-related life cycle assessments. The average supply chain emissions are as high as 25\% of all emissions from mining required for battery production, as provided in Supplementary Fig. S4. The underestimation of the supply chain's contribution to pre-operation emissions could be attributed to the inherent properties of supply chain analysis, which, in its nature, is complex. The analysis requires tracking of supply links, phases, modes, routes, product assembly details, and several other factors, making the process extensive and often perceived as insignificant in the overall magnitude of life cycle emissions. However, EV supply chain is quite particular due to the intensity of minerals requirements and their vast geographical distribution. Moreover, as a relatively new industry, with technology and resources concentrated in specific countries, the EV supply chain is therefore prone to higher emissions.

We can observe the growing complexity of the supply chain as we monitor the PMFs evolution from unimodal to multimodal PMF as the supply chain develops across phases. This results from the combinatorial impact of the transition between phases and the increased number of actors at each stage. However, PMFs transition back to unimodal PMF in the last phase of the supply chain. This reversion reflects the minimal international trade in the last phase, mainly due to high domestic demands, limited production capacities and EV tariffs on certain origins to protect local automobile industries.  

As for the emissions breakdown, it is clear that the last phase of the supply chain, namely the VM phase, is responsible for most of the EV supply chain's emissions. This can be attributed to 1) the high weight per KWh for vehicle transport, increasing the overall mass flow requirements per KWh compared to previous stages. 2) Emission factor per unit travel per weight. The emission factor for vehicle transport has a considerable higher unit value compared to other phases. This indeed is related to the type of transport vessel used for vehicles transport, which has a lower density-to-weight ratio when compared to the previous transported products at earlier supply phases. 

The breakdown also details the contribution of the transport mode per supply phase. A clear trend can be observed, with land transport having the highest percentage for mining and refining stages before sea transport takes over for the remainder of the supply chain. At the mining stage, we understand that this is mainly due to the location of the mines, which is beyond control and, in many cases, located in remote areas with large distances to cover to the nearest port for shipping. However, at the refining stage, this could be attributed to refining minerals in the producing countries or a result of the Chinese dominance over refining and battery production hubs, where most of the transport is then inland between these two phases. For the remainder of the supply chain, more actors require international shipping from the battery production hubs to the vehicle production hubs and thereafter to markets, where larger distances are covered, emitting higher emissions compared to the land distance, which is relatively low and closer to shipping ports. 

As part of the emissions breakdown, it is worth noting that the estimated emission values per battery and vehicle manufacturer do not reflect the company's production practices but rather their facility's geographical distribution with reference to the resources and markets. An important observation is that Hyundai's average supply chain emissions are considerably higher than those of other car manufacturers. This perhaps can be explained by considering the locations of Hyundai's production hubs compared to its share of different markets. For instance, Hyundai's home; South Korea is the only major EV actor with export values to the EU and US higher than its domestic inflow of KWh. This is coupled with the company's 10-11\% higher average vehicle weight per KWh compared to its counterparts, Tesla and VW, for example. The average vehicle weight for Hyundai is comparable to that of BYD, but BYD can maintain lower average supply chain emissions due to its localized market and supply chain. The impact of the localized supply chain on the emission values becomes more apparent when considering the emissions per market, where the Chinese market has significantly lower emissions than other markets, reflecting the strength of the localized Chinese supply chain and its positive impact on the emission values.    

The benefit of localization and regional supply structure in reducing supply chain emissions is a critical environmental consideration and a necessity to improve supply chain resilience. The resilience analysis is very clear and unequivocally demonstrates the Chinese dominance over the supply chain at all phases. The positive mass balance for China, coupled with its low percentage of non-domestic flow, indicates the high resilience of its supply chain, where the 20\% of non-domestic flow is mainly in the form of exports. In contrast, North American and EU's negative mass balances and high percentage of non-domestic flows (33\% and 40\%) indicate its high vulnerability to international imports with limited capacities to locally meet its demand at different supply levels. On the other hand, China took a very strategic position in the EV market, benefiting from the global shift toward electrification and the emergence of new technology to penetrate the automobile industry. This was not possible for ICEV vehicles with the generational trust built by older global companies. China also invested significantly in sourcing locally and intentionally with long-term contracts to secure the raw minerals required for EV production. In addition, China developed high-quality technology for mineral refining and, later on, reliable, high-quality battery and vehicle production lines. Some might argue that the Tesla giga factory in Shanghai contributed significantly to the improved quality of China's vehicle and battery production lines, as suppliers upped their quality standards to qualify for Tesla bids, and as a result, other local Chinese manufacturers emerged benefiting from high-quality supplies. 

Although this structure might be in favour of certain regions, the current structure puts the electrification strategy and the automobile industry at risk in other major markets. Hence, many initiatives exist, such as the Inflation Reduction Act (IRA) in the USA, the European Battery Alliance in the EU, and the Faraday Challenge in the UK. These initiatives will not only improve the supply chain's resilience but also significantly improve the environmental sustainability of the electrification strategy.


Additionally, the optimization exercise show promising potential for significant reductions in EV supply chain emissions. These reductions result from better distribution of the production hubs and the resource association that reduces the overall mass flow and emissions. Two scenarios were compared with the optimized market structure, one following the current global market structure and the other anticipating slight changes in international trade in the near future with increased export levels, especially for Chinese vehicles. Despite heavy tariffs, Chinese vehicles are still a very competitive option in a global market where they are, in many cases, a third of their competitor's price. In addition, China is the only country operating well below its production capacity at the moment, in anticipation of its expansion to regional and global markets in the future. However, the reduction based on the optimized supply chain for different market segments is significant for both current and near-future anticipated market dynamics.  

The reduction potential for the European market is the highest, mainly due to the geographical proximity of the market compared to the others. The supply chain's emissions for the Americas remained the highest after optimization, although it had a reduction percentage higher than the Asian market. The Asian market had the lowest, yet 32\% reduction potential. This is mainly because it is the closest to optimal structure with the advanced localized network. 
Overall, with the increase in demand and global efforts to speed up EV adoptions, considering optimal or near-optimal hub allocation and sourcing offers a serious potential not only for emissions reductions but also for faster, shorter, more efficient, and resilient supply chain architecture.

\section*{Methods}\label{sec4}

In this section, we provide a detailed breakdown of the mathematical formulation underlying the programmatic simulation model as well as the data structures used in the design and analysis. Additionally, this section provides details on the P-hub problem formulation for the supply chain optimization and resource allocation. These three dimensions are presented in three consecutive sub-sections hereunder.

\subsection*{Simulation Model Formulation}

The Transport emissions associated with the movement of goods can generally be expressed as a function of distance, weight and emission factor, as shown in Eq. \ref{eq1}. This simplistic representation acts as the framework for modelling the primary factors driving the emission values in the logistics chain, after which it's further refined to achieve a more accurate estimation of the supply chain environmental impacts incorporating available operational data.\\

Let \( E \) represent the emissions from the logistics supply chain, then \( E \) in its simplest form can be quantified by the function:

\begin{equation}
E = f(d, w, g) \label{eq1}
\end{equation}
where, 

\begin{itemize}
    \item \( d \): Distance travelled by the cargo between different supply phases. 
    \item \( w \): Weight of the cargo, which differs at each supply phase as it evolves towards the complete vehicle. 
    \item \( g \): Emission factor, which quantifies the amount of emissions produced per unit of weight per unit of distance, adjusted for the specific transportation mode.

\end{itemize}

The production of the EV has five primary phases with four mass flow operations. Let $(x)$ represent the supply phase, where $(x) \in \{E, P, B, V, M\}$, and where: 

\begin{itemize}
    \item \( E\): Extraction phase, which refers to minerals extraction locations. 
    \item \( P \): Processing phase, which refers to minerals processing locations.
    \item \( B \): Battery production, which refers to battery production locations. 
    \item \( V \): Vehicle production, which refers to Vehicle production locations. 
    \item \( M \): Market, which refers to the Market or the final destination of the produced vehicle. 
\end{itemize}

The location choice for phase \( E\), and \( P\) have \( k\) independent decisions, i.e., one decision per mineral, where \( k\) = number of minerals. Let $I_{x,i}$ be the index of the selected choice in phase $x$ for the $i$-th decision, which can be obtained from Eq. \ref{eq2}  for \( x \) $\in$ \{E, P, B\} \text { so that: } 

\begin{equation}
    I_{x,i} = \min \left\{ j : U_{x,i} \leq F_{j,x} \right\} \label{eq2}
\end{equation}

where \( U_{x,i} \sim \text{Uniform}(0, 1) \) is a random number, and \( F_{j,x} \) Cumulative probability distribution function for the $j$-th option for the $i$-th decision at phase $x$, and $i \in$ \{1,2,..,K\} \text{ for } \(x\) $\in$ \{E, P\} \text{ , $i$ = \( k\) = 1 for } \( x \) $\in$ \{B\}. \\

When $x \in \{V, M\}$, { Eq. \ref{eq3} with conditional choice set between two adjacent phases applies:}   

\begin{equation}
I_{x,i} = \min \left\{ j \in \mathcal{X}_{x}(I_{x-1,i}) : U_{x,i} \leq F_{j,x} \right\} \label{eq3}
\end{equation}
where:
\begin{itemize}
    \item \( I_{x-1,i} \) The index of the selected choice at the previous phase \( x-1 \), which influences the choice set \(\mathcal{X}_{x}\) of available choices at the current phase \( x \).
    \item \( U_{x,i} \) is a uniformly distributed random number, determining the probabilistic selection at phase \( x \).
    \item \( F_{j,x} \) is the cumulative probability function for the \( j \)-th option at phase \( x \), dependent on the set of options made available by the outcome of phase \( x-1 \).
    \item \( \mathcal{X}_{x}(I_{x-1,i}) \) represents the set of choices available at phase \( x \), which depends on the selected choice \( I_{x-1,i} \) from the previous phase.
\end{itemize}

Each mass transferred (\(w_{E_i}\)) between the selected choice $( I_{x,i} )$ at supply phase \( x \) and the next selected choice  $( I_{y,i} )$ at supply phase  \( y \), where \( x \) and \( y \) $\in$ \{E, P, B, V, M\} has two transport components; sea and land transport components as expressed in Eq. \ref{eq4} with associated emission factors per unit distance and unit weight that is adjusted to the transport mode and the transport vessel for each transport link, as shown in Eq. \ref{eq5} and \ref{eq6}: 

\begin{equation}
d_{I_{x,i}, I_{y,i}}^{total} = d_{I_{x,i}, I_{y,i}}^L + d_{I_{x,i}, I_{y,i}}^S
\label{eq4}
\end{equation}

\begin{equation}
g(T_{I_{X,i}, I_{Y,i}}^S) = \begin{cases} 
\gamma_1 & \text{if Bulk Carrier is used for sea transport between } I_{X,i} \text{ and } I_{Y,i} \\
\gamma_2 & \text{if Container Ship is used for sea transport between } I_{X,i} \text{ and } I_{Y,i} \\
\gamma_3 & \text{if Vehicle Transport ship is used for sea transport between } I_{X,i} \text{ and } I_{Y,i}
\end{cases} 
\label{eq5}
\end{equation}

\begin{equation}
g(T_{I_{X,i}, I_{Y,i}}^L) = \begin{cases} 
\beta_1 & \text{if Heavy Good Vehicle / Diesel is used for land transport between } I_{X,i} \text{ and } I_{Y,i} \\
\beta_2 & \text{if Articulated - Vehicle Transport is used for land transport between } I_{X,i} \text{ and } I_{Y,i}
\end{cases} 
\label{eq6}
\end{equation}

Incorporating the details from Eq. \ref{eq4}, \ref{eq5} and \ref{eq6} into Eq. \ref{eq1}, the emission value $E$ can then be obtained using Eq. \ref{eq7} incorporating distances, weights, and emission factors for transitions between all supply phases of the logistic chain with $k$+3 decision choices made for each emission value calculated:

The emissions equation is formulated as follows:

\begin{align}
E = & \sum_{i=1}^K \Bigg( w_{E_i} \cdot \left(d_{I_{E,i}, I_{P,i}}^L \cdot g(T_{I_{E,i}, I_{P,i}}^L) + d_{I_{E,i}, I_{P,i}}^S \cdot g(T_{I_{E,i}, I_{P,i}}^S)\right) \nonumber \\
& + w_{P_i} \cdot \left(d_{I_{P,i}, I_{B}}^L \cdot g(T_{I_{P,i}, I_{B}}^L) + d_{I_{P,i}, I_{B}}^S \cdot g(T_{I_{P,i}, I_{B}}^S)\right) \nonumber \\
& + w_B \cdot \left(d_{I_{B}, I_{V}}^L \cdot g(T_{I_{B}, I_{V}}^L) + d_{I_{B}, I_{V}}^S \cdot g(T_{I_{B}, I_{V}}^S)\right) \nonumber \\
& + w_V \cdot \left(d_{I_{V}, I_{M}}^L \cdot g(T_{I_{V}, I_{M}}^L) + d_{I_{V}, I_{M}}^S \cdot g(T_{I_{V}, I_{M}}^S)\right). \label{eq7}
\end{align}

\textbf{Where:}
\begin{itemize}
    \item $w_{X_i}$: Weight of the cargo at phase $X$ for selection $i$.
    \item $d_{I_{X,i}, I_{Y,i}}^L$: Land distance between the selections $i$ from phase $X$ to phase $Y$.
    \item $d_{I_{X,i}, I_{Y,i}}^S$: Sea distance between the selections $i$ from phase $X$ to phase $Y$.
    \item $g(T_{I_{X,i}, I_{Y,i}}^L)$: Emission factor for land transport between the selections $i$ from phase $X$ to phase $Y$.
    \item $g(T_{I_{X,i}, I_{Y,i}}^S)$: Emission factor for sea transport between the selections $i$ from phase $X$ to phase $Y$.
    \item $I_{X,i}$: Index of the selected choice at phase $X$ for the $i$-th decision.
    \item $U_{X,i}$: Uniformly distributed random number for phase $X$ used to select $I_{X,i}$.
    \item $F_{j,X}$: Cumulative probability distribution function for the $j$-th option for the $i$-th decision at phase $X$.
\end{itemize}

Equation \ref{eq7} calculates the emission value for each interaction, after which the cumulative average is derived over $N$ iterations, as shown in Equation \ref{eq8}, considering only a subset of choices $k$ that account for varying mineral compositions across different battery types. Additional control parameters are applied in the programmatic model, especially for the analysis of the simulation results to associate different battery and vehicle manufacturers with specific trade links, cargo weights and battery compositions. 

\begin{equation}
\bar{E} = \frac{1}{N} \sum_{n=1}^N E_{n, \mathcal{S}} \label{eq8}
\end{equation}

where:
\begin{itemize}
    \item \( E_{n, \mathcal{S}} \) represents the emissions from the \( n \)-th iteration, including only the emissions corresponding to the specific subset \( \mathcal{S} \subset K \) of choices according to the battery type.
    \item \( N \) is the total number of iterations over which the average is calculated.
    \item \( \mathcal{S} \) is defined within the set \( K \), such as including the choices that are part of a given battery type composition.  
\end{itemize}

\subsection*{ Simulation Model Inputs} 

\begin{figure}[ht!]
  \centering
  \includegraphics[height=9cm]{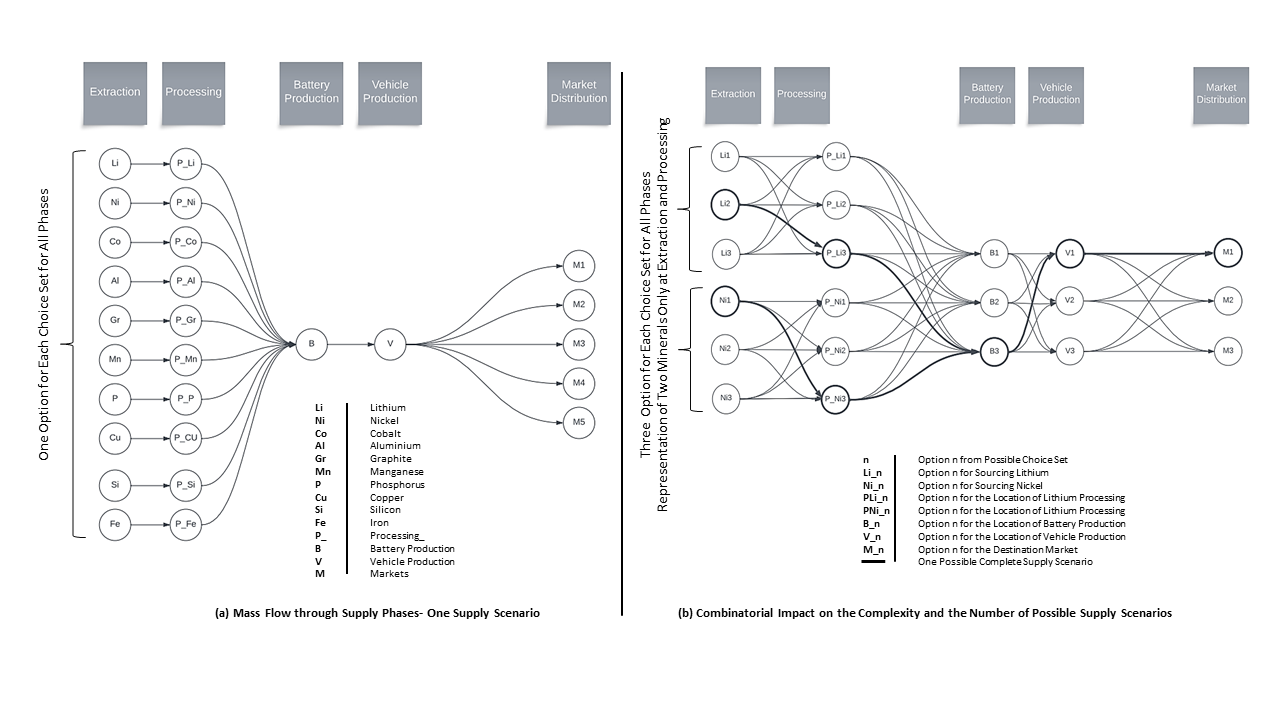}
  \caption{EV Production Mass Flow Combinatorial Scenarios Impact Illustration \textbf{(a)} Mass Flow through Supply Phases- One Supply Scenario \textbf{ (b)} Combinatorial Impact on the Complexity and the Number of Possible Supply Scenarios}
  \label{fig:6}
\end{figure}

The most significant challenge in evaluating supply chain emissions is, without a doubt, the access to reliable data. Origin and Destination (O$-$D) data for different phases of the logistics chain are essential to obtain the mass flow, and therefore the emission values. Specific O-D and transport mode data are very hard to obtain, and while aggregate-level data might be available, the competitive nature of the industry makes it very difficult to access detailed operational-level data. Alternatively, we could consider possible supply scenarios when exact supply links are not available. However, this approach becomes complicated and complex very quickly. For example, Fig. \ref{fig:6} (a) illustrates the logistics supply process required for the production of EVs, showcasing one possible supply scenario with a single source option for each part of every phase of the supply chain. The problem arises due to the large number of minerals required and the geographical distribution of these mineral origins for battery and vehicle production. This interrelation of different suppliers at various stages results in a vast number of possible supply scenarios, especially when conducting a global-scale evaluation. To illustrate this, Fig. \ref{fig:6} (b) shows the different flow scenarios only for lithium and nickel, with three possible choices at each level, demonstrating how quickly the combinatorial alternatives at each stage increase the complexity of the EV supply chain. The bold lines represent just one supply scenario for both lithium and nickel out of $3^7$ possible scenarios depicted in the graph.    

In this research, we used industry reports, company profiles, governmental and intergovernmental reports, and data from consulting firms to establish the O-D matrix for each supply phase. The total number of possible supply scenarios based on the constructed O-D matrix amounts to 2.14 trillion scenarios. To tackle this complexity, we developed the SCEV-Sim model to identify and simulate the most probable supply scenarios from the set and calculate the supply emissions for different markets and manufacturers. Table \ref{tab: SCEV Model Data Inputs Summary} outlines the sources for the most relevant inputs used in constructing the O-D matrix and details SCEV-Sim's market coverage across each supply phase.  

\begin{table}[H]
\caption{SCEV-Sim Data Inputs Summary}
\centering
\renewcommand{\arraystretch}{0.5}
\footnotesize
\begin{tabular*}{\textwidth}{@{\extracolsep{\fill}}|p{0.2\textwidth}|p{0.35\textwidth}|p{0.1\textwidth}|p{0.1\textwidth}|p{0.075\textwidth}|@{}}
\hline
\textbf{Data Type} &  \textbf{Description} &  \textbf{Related Supply Phase} &  \textbf{Overall Phase \% Market Coverage}  & \textbf{Source} \\
\hline
&&&&\\
Government Federal Reports & US geological commodity summaries, breakdown of global mineral production values per country, mineral reserves, main trends, prices, etc. &  Extraction &  90\%& \cite{usgs2023minerals} \\
&&&&\\

Data analytics and Consulting Firm Reports & GlobalData mining reports on commodity outlooks, mining projects and their locations and productions details and global share &  Extraction &  90\% & \cite{GlobalBauxiteMining2023, GlobalCobaltMining2022, GlobalCopperMining2022, GlobalMetalsMining2022, GlobalIronOreMining2022, GlobalLithiumMining2022, GlobalGraphiteMining2023} \\
&&&&\\
Government Federal Report& US National Geospatial-Intelligence Agency (NGA) publication, detailed ports database with world port index, port infrastructure, handling details, and location coordinates & Interphase Transfer & NA & \cite{UpdatedPub150CSV}\\
&&&&\\       

Software Library & Maritime Sea route distances. A Python package to generate the shortest maritime sea route for sailing between two ports & Interphase Transfer & NA & \cite{searoute}\\
&&&&\\

Inter-governmental Agency Reports & IEA intergovernmental agency analysis based on US geological commodity summaries, S\&P global and Benchmark Mineral Intelligence and Wood Mackenzie &  Processing & Li 99\%, Co 88\%, C 100\%, Ni 65\%, Cu 57\% & \cite{iea2023processing}\\
&&&&\\
Inter-governmental reports, Company profiles, industry reports & Intergovernmental Agency reports on annual market shares. Company profiles reveal production facilities' geographical distribution. Industry reports on annual production volumes & Battery production hubs & 79.6\% & \cite{CATLProfile, LGEnergyGlobalNetwork, PanasonicGlobalNetwork,EVVolumes2023, byd2024globalization, iea2024global}\\ 
&&&&\\
Inter-governmental reports, Company profiles, industry reports & Intergovernmental Agency reports on annual market shares. Company profiles reveal production facilities' geographical distribution. Industry reports on annual production volumes & Vehicle production hubs& 60\% & \cite{iea2024global, tesla2024, hyundai2024, byd2024globalization, evvolumes2023sales} \\ 
&&&&\\
Inter-governmental reports, Company profiles, industry and consulting company reports & Intergovernmental Agency reports on mass flows. Company profiles reveal sales per region, consulting company reports on export and import origin-destination and volumes& Interphase Mass transfer joint and conditional probabilities & NA & \cite{iea2024global, ICCT_MajorMarkets2023,OEC_ElectricMotorVehicles, HyundaiSalesResults, Volkswagen2023Growth} \\
\hline
\end{tabular*}
\label{tab: SCEV Model Data Inputs Summary}
\end{table}

\subsection*{Global Supply Chain Optimization and Resource Allocation}\

To enhance the sustainability and resilience of the EV supply chain and ensure that car manufacturers have viable options for production hubs that contribute to a more sustainable and resilient supply structure, a P-Hub optimization problem was used to identify potential hub locations across different regions, while optimizing the hubs sourcing options for minerals required in the production process.

The formulation of the optimization problem begins with defining the optimization space. Since the flow moves from sources to manufacturing hubs and then to markets, we can define the optimization space as follows:

\begin{itemize}
    \item \( J = \{1, \dots, j_{max}\} \): the set of subsets of sources and $j_{max}$ is the total number of different minerals that compose the three battery types. 
    \item \( I = \{1, \dots, \alpha\} \): the set of potential hubs.
    \item \( M = \{1, \dots, \beta\} \): the set of markets.
    \item \( J_j = \{1, \dots, n_j\} \): the set of options in subset $j$, where $n_j$ varies by subset. The set of options in subset $j$ represents the number of sources available for each mineral. 
\end{itemize}

It is important to note here, as we consider the optimal supply chain structure, two key assumptions are made: 

\begin{enumerate}
    \item Processing of minerals is not a standalone hub, where it either takes place at the same location as extraction or battery production hubs. This is only valid for an optimal market structure and for minimizing the overall travel distance, cost and emissions. This assumption only applies to the optimization of the future market, while it does not apply to the analysis of the current market structure. It is thus not used in the analysis section or the simulation model. 
    \item There is one level of production hubs. In other words, batteries and vehicles are produced in the same hub or with no significant distance between the two production hubs. This structure is assumed to achieve global optimum scenarios for each market, and while it is only assumed for the optimization phase, it is already adopted in cases like Canada's localization approach and BYD's in-house battery and vehicle production.  
\end{enumerate}

For the optimization of the supply chain, two cost matrices are constructed for each market scenario. The first matrix $c_{ijk}$ comprises the cost of transporting the minerals from the extraction location to the production hub. The second matrix $c_{km}$ provides the cost of transferring the vehicles from the production hubs to the prospective markets as detailed in the following description of the variables: 

\begin{itemize}
    \item \( c_{ijk} \): cost of transporting minerals from option \( i \) in subset \( j \) to hub \( k \).
    \item \( c_{km} \): cost of transporting electric vehicles from hub \( k \) to market \( m \).
\end{itemize}

Two binary decision variables are introduced into the optimization space: one for sourcing options $x_{ijk}$ and another for the production hubs $y_k$. These variables indicate whether a specific hub is active and whether a specific transportation link (from a source option to a hub) is utilized to achieve an optimal supply chain structure as follows. 

\begin{itemize}
    \item \( x_{ijk} \in \{0,1\} \): binary variable that equals 1 if option \( i \) from subset \( j \) is shipped through hub \( k \), and 0 otherwise.
    \item \( y_k \in \{0,1\} \): binary variable that equals 1 if hub \( k \) is selected, and 0 otherwise.
\end{itemize}

With the objective of minimizing the overall transportation cost (emissions), the objective function is set in Eq. \ref{eq9} as follows: 

\begin{equation}
    \text{Minimize} \quad Z = \left( \sum_{k=1}^{\alpha} \left( \sum_{j=1}^{j_{max}} \sum_{i=1}^{n_j} c_{ijk} x_{ijk} \right) \right) + \left( \sum_{k=1}^{\alpha} \sum_{m=1}^{\beta} c_{km} y_k \right) \label{eq9}
\end{equation} \

Where:
\begin{itemize}
    \item \( c_{ijk} \) is the transportation cost from option \( i \) of subset \( j \) to hub \( k \).
    \item \( c_{km} \) is the transportation cost from hub \( k \) to market \( m \).
    \item \( x_{ijk} \) indicates whether option \( i \) of subset \( j \) is transported through hub \( k \).
    \item \( y_k \) indicates whether hub \( k \) is selected for use.
\end{itemize}

Two constraints are added to the objective function. The first constraint ensures that exactly two hubs are selected for each market. This constraint is represented in Eq. \ref{eq10}: 
\begin{equation}
    \sum_{k=1}^{\alpha} y_k = 2
    \label{eq10}
\end{equation} 

The second constraint ensures that one option from each \( j \) is selected for each hub \( k \). This is meant to ensure that at least one sourcing option is selected for each of the minerals required for the production of the EV. This constraint is represented in Eq. \ref{eq11}:  

\begin{equation}
    \sum_{i=1}^{n_j} x_{ijk} = y_k \quad \forall j, \forall k
    \label{eq11}
\end{equation}

Finally, Eq. \ref{eq12} links the selected hub with the sourcing option, ensuring that option \( i \) from subset \( j \) can be selected for transportation through the hub \( k \) only if hub \( k \) is selected.
\begin{equation}
    x_{ijk} \leq y_k \quad \forall i, \forall j, \forall k \label{eq12}
\end{equation}

This mixed integer linear programming model (MILP) established in Eqs. \ref{eq9} to \ref{eq12} is addressed using the coin-or-branch and cut (CBC) solver. The CBC's branch-and-cut algorithm was used due to its high effectiveness in managing complex problems and suitability for large-scale MILP since the algorithm is enhanced using cutting plane techniques and advanced heuristics that guide the search and allow for rapid convergence to optimal solution despite the extensive combinatorial elements. 

\section*{Data Availability}

\begin{itemize}
\item \textbf{Data availability}:  Available on the project GitHub repository: \url{https://github.com/LiTrans/SCEV-Sim}  Access for reviewers to the full dataset is provided via supplementary files. \\
\item \textbf{Code availability:} Available on the project GitHub repository \url{https://github.com/LiTrans/SCEV-Sim}  Access for reviewers to the full code scripts is provided via supplementary files. \\
\end{itemize}
\bibliography{sample}

\section*{Acknowledgements}

\textbf{Funding} This research was funded by a grant from the Canada Research Chair program in Disruptive Transportation Technologies and Services (CRC-2021-00480) and NSERC Discovery (RGPIN-2020-04492) fund.\\

\section*{Author contributions statement}

The authors confirm their contribution to the paper as follows: \textbf{Tareq Alsaleh}.: Conceptualization, Methodology, Data curation, Investigation, Formal analysis, Software, Visualization, Writing - original draft, and Writing - review \& editing. \textbf{Bilal Farooq}.: Conceptualization, Methodology, Investigation, Funding acquisition, Project administration, Resources, Supervision, and Writing - review \& editing. 

\section*{Additional information}

\textbf{Competing interests:} The authors declare no competing interests.

\end{document}